\documentclass[12pt]{article}
\usepackage{amsthm, amsmath}
\usepackage{fullpage}
\usepackage{xcolor}
\usepackage[colorlinks=false,linkbordercolor=red]{hyperref}
\hypersetup{%
  colorlinks=false,
  linkbordercolor=red,
  pdfborderstyle={/S/U/W 1}
}

\def\WE{\widehat{E}}
\def\SUM{t}
\def\l{\ell}
\def\GLi{G[V_i,V_{i-1}]}
\def\GsLi{G_{\sigma}[V_i,V_{i-1}]}
\def\MT{\mathcal{T}}

\newtheorem{conj}{Conjecture}
\newtheorem*{theoremA}{Main Theorem}
\newtheorem*{keylemma}{Helpful Lemma}

\begin{document}
\author{Daniel W. Cranston\thanks{Department of Mathematics and Applied Mathematics, Virginia Commonwealth University, Richmond, VA, 23284. email: \texttt{dcranston@vcu.edu}} 
}
 
\title{Regular graphs of odd degree are antimagic}
\maketitle
\begin{abstract}
An antimagic labeling of a graph $G$ with $m$ edges is a bijection
from $E(G)$ to $\{1,2,\ldots,m\}$ such that for all vertices $u$ and $v$, the
sum of labels on edges incident to $u$ differs from that for edges incident to
$v$.  Hartsfield and Ringel conjectured that every connected graph other than
the single edge $K_2$ has an antimagic labeling.  We prove this conjecture for
regular graphs of odd degree.
\end{abstract}
\section{Introduction}
A \emph{magic square} of order $n$ is a $n\times n$ arrangement of the integers
$\{1,2,\ldots, n^2\}$ so that the sums of the entries in each row, each column,
and along the two main diagonals are equal.  These squares were known to the Chinese as
early as the fourth century B.C.\ and have been widely studied in recreational
mathematics~\cite{Gar88}.

A \emph{labeling} of a graph $G$ with $m$ edges is a bijection from $E(G)$ to
$\{1,2,\ldots,m\}$.  Given a labeling of a graph, the \emph{vertex sum} at a
vertex $v$ is the sum of the labels on edges incident to $v$.
A labeling is \emph{magic} if all vertex sums are equal.  Magic labelings take
their name from their connection with magic squares, since a magic square of order
$n$ naturally gives rise to a magic labeling of the complete bipartite graph
$K_{n,n}$ (vertices in one part correspond to rows of the square, and 
vertices in the other correspond to columns).  Finally, a labeling of a graph is
\emph{antimagic} if all its vertex sums are different.  We call a graph
antimagic (magic) if it has an antimagic (magic) labeling.

It is easy to find many graphs that are not magic (for example, forests). 
However, graphs that are not antimagic are rare.  In fact, Hartsfield and Ringel 
conjectured the following.

\begin{conj}[\cite{HR90}]
Every connected graph other than $K_2$ is antimagic.
\end{conj}

Hartsfield and Ringel also explicitly conjectured that all trees other than
$K_2$ are antimagic.  Both conjectures remain wide open; however, much progress
has been made.  The first major result on antimagic labelings was due to 
Alon, Kaplan, Lev, Roditty, and Yuster~\cite{AKLRY04}. 
They showed that there exists a constant $c$ such that if $G$ is an $n$-vertex
graph with minimum degree $\delta \ge c\log n$, then $G$ is antimagic.
This proof relies on a combination of combinatorial ideas, probabilistic
tools, and methods from analytic number theory.
They also proved that graphs with maximum degree $\Delta\ge n-2$ are antimagic.
Yilma~\cite{Yil13+} later extended this result to show that graphs with  $\Delta\ge n-3$ 
are antimagic.  His proof finds a breadth-first spanning tree $T$ rooted at a vertex
of maximum degree; he labels all edges outside of $T$ first, then uses the
largest $n-1$ labels on $T$ to guarantee an antimagic labeling.

Hefetz~\cite{Hef05} used algebraic tools to show that a graph is antimagic if it has $3^k$
vertices and a $C_3$-factor.  Hefetz, Saluz, and Tran~\cite{HST10} generalized
this approach to show that a graph is antimagic if it has $p^k$ vertices and a
$C_p$-factor (where $p$ is an odd prime).  Cranston~\cite{Cra09} used Hall's marriage theorem
to show that regular bipartite graphs are antimagic.  
Liang and Zhu~\cite{LZ13+} labeled edges in order of decreasing distance from a
central vertex (breaking ties carefully) to show that 3-regular graphs are antimagic.

Perhaps the most interesting result is that of Eccles~\cite{Ecc13+}, who
recently improved on the work of Alon et.al.  He showed that if a graph has no
isolated edges or vertices and has average degree at least
4468, then it is antimagic.  
He conjectures that, under the same first condition,
average degree at least $\sqrt{2}$ implies that a graph is antimagic.
This much stronger conjecture immediately implies Conjecture~1, since a
connected $n$-vertex graph has at least $n-1$ edges, and so for $n\ge 4$
has average degree at least $2(n-1)/n=2-2/n>\sqrt{2}$. 

In this note, we prove that every $k$-regular graph with $k$ odd and $k\ge 3$ is
antimagic.

\section{Main Result}
A \emph{trail} is a walk in $G$ that may reuse vertices but may not reuse edges;
a trail is \emph{open} if it starts and ends at distinct vertices, and is  
\emph{even (odd)} if its length is even (odd).
For a set of vertices $U$ and a function $\sigma$ we write $\sigma(U)$ to denote
$\{\sigma(u):u\in U\}$.  For a subgraph or trail $H$, we write $d_H(v)$ for the
degree of $v$ in $H$.  
We begin with an easy decomposition result for bipartite graphs.

\begin{keylemma}
\label{lemma1}
Let $G$ be a bipartite graph with parts $U$ and $W$.  There exists a
function $\sigma: U\to E(G)$ and a set $\MT=\{T_1,T_2,\ldots\}$ such
that $\sigma(u)$ is incident to $u$ for all $u\in U$ and $\MT$ is a
collection of edge-disjoint open trails with at most one trail ending at each
vertex and with $(\bigcup_{T\in\MT}E(T))\cap \sigma(U)=\emptyset$ and
$\bigcup_{T\in\MT}E(T)\cup \sigma(U) = E(G)$.  In other words, we can
partition $E(G)$ into $\MT$ and $\sigma(U)$.
\end{keylemma}
\begin{proof}
We first choose $\sigma(U)$ arbitrarily, and let $\WE=E(G)\setminus
\sigma(U)$.  We form a greedy trail decomposition of $\WE$ as follows.  Start at an
arbitrary vertex and keep walking (using unused edges of $\WE$) as long as possible.
When you reach a vertex with no unused edges, start a trail at another vertex. 
Repeat this process until all edges are used up.  This gives a decomposition
$\MT$ of $\WE$, but it might contain a closed trail.

Suppose that $\MT$ contains a closed trail $T_1$. If any
vertex $v$ of $T_1$ has an open trail $T_2$ that ends at $v$, then we 
splice $T_1$ and $T_2$ together, by starting at $v$, following all the edges of
$T_1$, then following the edges of $T_2$.  
If no vertex of $T_1$ is the endpoint of an open trail in $\MT$,
then choose $u\in U\cap V(T_1)$ arbitrarily.  Let $w$ be a successor of $u$ on
$T_1$ and let $v$ be such that $\sigma(u)=uv$.  We redefine $\sigma(u):=uw$,
and redefine $T_1:=T_1-uw+uv$.
Now $T_1$ is an open trail, since $d_{T_1}(w)$ is odd.  

By repeating this process for each closed trail in $\MT$, we reach a collection
of open trails.  If any vertex $v$ is the endpoint of at least two open trails,
then we merge them together, by walking along one to end at $v$, then walking
along another starting from $v$.  Merging two trails reduces the number of
trails, and ``opening up'' a closed trail (as described above), does not
increase this number.  So iterating these merging and opening up steps
gives the desired partition of $E(G)$ into $\MT$ and $\sigma(U)$.
%
%
%
%
\end{proof}

Now we prove our main result.
Our proof builds heavily on that of Liang and Zhu~\cite{LZ13+}, who showed that
3-regular graphs are antimagic.  

\begin{theoremA}
Every $k$-regular graph with $k$ odd and $k\ge 3$ is antimagic.
\end{theoremA}
\begin{proof}
Suppose that $G$ and $H$ are both antimagic $k$-regular graphs and that
$|E(G)|=m$. 
Given antimagic labelings for $G$ and $H$, we get an antimagic labeling
of $G\cup H$ by increasing the label on each edge of $H$ by $m$.  
Thus, we need only consider connected graphs.

Choose an arbitrary vertex $v^*$ and let $V_i$ denote the set of vertices at
distance exactly $i$ from $v^*$; let $p$ be the furthest distance of a vertex
from $v^*$.  Let $G[V_i]$ denote the subgraph induced by $V_i$ and 
$\GLi$ denote the induced bipartite subgraph with parts
$V_i$ and $V_{i-1}$.

For each $i$, we apply the \hyperref[lemma1]{Helpful Lemma} to $\GLi$ with $U=V_i$
and $W=V_{i-1}$ to get a partition of $E(\GLi)$ into an edge set
$\sigma(V_i)$ and a collection of edge-disjoint open trails.  Let
$\GsLi = \GLi\setminus \sigma(V_i)$.  Let
$E_i=E(G[V_i])$, let $E'_i=E(\GsLi)$, and let
$E''_i=\sigma(V_i)$; note that $E'_i$ and $E''_i$ partition $E(\GLi)$.
Given a labeling $f$ of the edges, we denote the total sum of labels on edges
incident to vertex $v$ by $\SUM(v)=\sum_{e\in E(v)}f(e)$, where $E(v)$ denotes
the set of edges incident to $v$.  Similarly, we denote the partial sum at $v$ 
(omitting the label on $\sigma(v)$) by
$p(v)=\sum_{e\in E(v)\setminus\{\sigma(v)\}}f(e)=\SUM(v)-f(\sigma(v))$.

We now outline the proof.  We will label the edges in the order $E_p, E'_p,
E''_p, \ldots, E_1, E'_1, E''_1$, using the smallest unused labels on each edge
set when we come to it.  In other words, we use the $|E_p|$ smallest labels on
$E_p$, the $|E'_p|$ next smallest labels on $E'_p$, the $|E''_p|$ next smallest
labels after that on $E''_p$, etc.  (Note that the labels assigned to each of
these edge sets span an interval.)  This label assignment immediately gives
that if $i\ge j+2$ and $u\in V_i$ and $w\in V_j$, then $\SUM(u)<\SUM(w)$ since
$G$ is regular and the edges incident to $u$ have smaller labels than the edges
incident to $w$.  Thus, we need only ensure that $\SUM(u)\ne \SUM(w)$ when
either (i) $u,w\in V_i$ or (ii) $u\in V_i$ and $w\in V_{i-1}$.  We handle these
two cases by specifying more precisely how to assign the label to each edge of
these $3p$ edge sets.  

We label the edges of each $E_i$ arbitrarily from its assigned labels.
We now specify how to label each $E''_i$; in the process, we handle Case (i).
Suppose that for some $i$, we have already labeled the edges of $E_p, E'_p,
E''_p, \ldots, E_i, E'_i$. As a result, $p(u)$ is already determined for each
$u\in V_i$.  We may name the vertices of $V_i$ as $u_1, u_2, u_3, \ldots$ so
that $p(u_1)\le p(u_2) \le p(u_3) \le \cdots$.  Now we use the smallest label
for $E''_i$ on $\sigma(u_1)$, the next smallest on $\sigma(u_2)$, etc.  This
ensures that $\SUM(u_j)<\SUM(u_{j+1})$ for all $u_j\in V_i$.

Finally, we specify how to label each $E'_i$; in the process, we handle Case
(ii).  That is, we ensure that if $u\in V_i$ and $w\in V_{i-1}$, then
$\SUM(u)\ne\SUM(w)$.  Let $\{s,s+1,\ldots,\l-1,\l\}$ be the set of labels to be
used on $E'_i$.  Recall that $G$ is $k$-regular for odd
$k\ge 3$, and let $t=(k-1)/2$.  We will ensure that $p(u)\le t(s+\l)$ and that
$p(w)\ge t(s+\l)$.  Now since $f(\sigma(u))<f(\sigma(w))$, we get that $t(u)<t(w)$.
The details follow.

Let $\MT$ be the set of open trails partitioning $E'_i$ (from the
\hyperref[lemma1]{Helpful Lemma}).  Again, let
$\{s,s+1,\ldots,\l-1,\l\}$ be the labels assigned to $E'_i$.  We label
each trail so that every pair of successive labels (on a trail) incident to a
vertex $u\in V_i$ has sum at most $s+\l$ 
and each pair of successive labels incident to a vertex $w\in V_{i-1}$
has sum at least $s+\l$. 
This ensures that $p(u)\le t(s+\l)$ and $p(w)\ge t(s+\l)$.

We first label each even trail, then label the odd trails, taken together in
pairs (possibly with a single odd trail last).
Suppose that we have already labeled some even number $2r$ of edges in the set
$E'_i$ and the remaining labels available for this edge set are $\{s+r,
s+r+1,\ldots, \l-r-1,\l-r\}$. We have three possibilities. (1) Suppose first
that $T\in \MT$ is an
even trail with both endpoints in $V_{i-1}$.  We assign the labels: $s+r,
\l-r,s+r+1,\l-r-1,\ldots$ successively along the trail.  Now every two successive
edges incident to $u\in V_i$ have sum $s+\l$ and every two successive edges
incident to $w\in V_{i-1}$ have sum $s+\l+1$.  (2) Suppose instead that $T\in
\MT$ is an even trail with both endpoints in $V_i$.  Now we assign
the labels: $\l-r, s+r, \l-r-1, s+r+1, \ldots$ successively along the trail. 
Now every two successive edges incident to $u\in V_i$ have sum $s+\l-1$ and
every two successive edges incident to $w\in V_{i-1}$ have sum $s+\l$.  (3)
Finally, suppose that $T_1, T_2\in \MT$ are odd trails with lengths
$2a+1$ and $2b+1$.  Beginning at a vertex in $V_i$, we label the edges of $T_1$
with $\l-r, s+r, \l-r-1, \ldots, s+r+a-1, \l-r-a$.  Here the successive pairs
of labels incident to $u\in V_i$ sum to $s+\l-1$ and the pairs incident to
$w\in V_{i-1}$ sum to $s+\l$.  Finally, beginning at a vertex in $V_{i-1}$,
we label the edges of $T_2$ with $s+r+a, \l-r-a-1, s+r+a+1, \ldots, s+r+a+b$.
Again the successive pairs incident to $u\in V_i$ sum to $s+\l-1$ and the
successive pairs incident to $w\in V_{i-1}$ sum to $s+\l$.
If we have a single odd trail left at the end, we treat it like a trail of
length $2a+1$ above.

All that remains is to verify that for $u\in V_i$ and $w\in V_{i-1}$ we
have $p(u)\le t(s+\l)$ and $p(w)\ge t(s+\l)$.  We consider $p(u)$, and the
analysis for $p(w)$ is nearly identical.  Recall that $d_G(u)=k=2t+1$.  If
$d_{E'_i}(u)=2t$, then the desired inequality holds, since each of the $t$
pairs of successive edges on trails through $u$ have label sum at most $s+\l$. 
If $u$ is the end of some trail $T$ in $\MT$, then let $e$ be the final edge of
$T$ incident to $u$; note that $f(e)\le \l$.  But now, we have $d_{E'_i}(u)$ is
odd, so $u$ has some incident edge (in fact, an odd number of them) in 
$E'_{i+1}\cup E''_{i+1}\cup E_i$; this edge has label less than $s$.  Thus, the
sum of this label and $f(e)$ is less than $s+\l$.  If $u$ has additional
incident edges in  $E'_{i+1}\cup E''_{i+1}\cup E_i$, then each edge has label
less than $s$; thus, each pair of these edges has label sum less than $s+\l$. 
So $p(u)\le t(s+\l)$, as desired.  For each $w\in V_{i-1}$, the analysis to
show that $p(w)\ge t(s+\l)$ is nearly identical to that above; the only
difference is that all edges incident to $w$ that are not in $E'_i$ are in
$E''_i\cup E_{i-1}\cup E'_{i-1}$, so each such edge has label larger than $\l$.
This completes the proof.
%
\end{proof}

We remark in closing that the proof easily translates to an efficient
(polynomial time) algorithm to find an antimagic labeling.  
We thank Mike Barrus for his careful reading of this manuscript and detailed
feedback.

\end{document}